\documentclass[final,3p,times]{elsarticle}
\RequirePackage{geometry}
\usepackage{amssymb}
\usepackage{lmodern}        % Added this package to enable italic font

\usepackage{caption}
\usepackage{amsmath}
\usepackage{array}
\usepackage{pifont}
\usepackage{subcaption}
\usepackage{multirow}
\usepackage{etoolbox}
\usepackage{natbib}
\usepackage{color}
\usepackage[dvipsnames]{xcolor}
\usepackage{ulem}
\apptocmd{\sloppy}{\hbadness 10000\relax}{}{}
\usepackage{enumitem}
\usepackage{algorithm}
\usepackage{algpseudocode}
\usepackage{textcomp}
\usepackage{amsmath}
\usepackage{hyperref}
\usepackage{hyperref}
\usepackage{color}
\usepackage{amsmath}
\usepackage{amsmath}
\usepackage{algpseudocode}
\usepackage{algpseudocode}
\usepackage{algpseudocode}
\usepackage{amsmath}

\newcommand{\fetwo}{FE$^\mathrm{2}$}

\begin{document}

\begin{frontmatter}

\title{Mixing Data-Driven and Physics-Based Constitutive Models using Uncertainty-Driven Phase Fields}

\author[label1]{J. Storm}
\author[label2]{W. Sun}
\author[label1]{I. B. C. M. Rocha}
\author[label1]{F. P. van der Meer}

 \affiliation[label1]{organization={Delft University of Technology, Department of Civil Engineering and Geosciences},
             addressline={P.O. Box 5048},
             city={Delft},
             postcode={2600GA},
             state={Zuid-Holland},
             country={The Netherlands}}

\affiliation[label2]{organization={Columbia University, Department of Civil Engineering and Engineering Mechanics},
         addressline={Mail Code: 4709},
         city={New York},
         postcode={NY 10027},
         state={New York},
         country={United States of America}}

\begin{abstract}
There is a high interest in accelerating multiscale models using data-driven surrogate modeling techniques.
Creating a large training dataset encompassing all relevant load scenarios is essential for a good surrogate, yet the computational cost of producing this data quickly becomes a limiting factor.
Commonly, a pre-trained surrogate is used throughout the computational domain.
Here, we introduce an alternative adaptive mixture approach that uses a fast probabilistic surrogate model as constitutive model when possible, but resorts back to the true high-fidelity model when necessary.
The surrogate is thus not required to be accurate for every possible load condition, enabling a significant reduction in the data collection time.
We achieve this by creating phases in the computational domain corresponding to the different models.
These phases evolve using a phase-field model driven by the surrogate uncertainty.
When the surrogate uncertainty becomes large, the phase-field model causes a local transition from the surrogate to the high-fidelity model, maintaining a highly accurate simulation.
We discuss the requirements of this approach to achieve accurate and numerically stable results and compare the phase-field model to a purely local approach that does not enforce spatial smoothness for the phase mixing.
Using a Gaussian Process surrogate for an elasto-plastic material, we demonstrate the potential of this mixture of models to accelerate multiscale simulations.
\end{abstract}


\end{frontmatter}

\section{Introduction}\label{sec1}
Advanced manufacturing techniques allow for the creation of materials with properties tailored to their use.
Multiscale modeling can accurately capture the behavior of these materials, yet its use is limited by the high computational cost of simulation.
The need for more computationally efficient models, combined with the progress in machine learning, has resulted in numerous data-driven models~\citep{J.GhaboussiJ.H.GarrettJr.1991, ghavamian2019accelerating, raissi2019physics, lu2021learning, vlassis2023geometric, fuhg2024review}.
Given sufficient data, these models often show an impressive ability to capture the behavior of complex history-dependent load paths.
They can therefore be used as surrogates for the representative volume element (micromodel) in a multiscale analysis.
However, obtaining the necessary training data from microscale simulations is very costly.
Because purely data-driven methods might perform poorly in extrapolation, data requirements necessary to ensure fidelity tend to be high.
In addition, it is often challenging to predict which loading scenarios will be relevant.
The surrogate therefore needs to be trained across a broad range of inputs to generalize for unseen data, many of which may never be relevant during the target simulation.

We can broadly categorize two (non-exclusive) approaches to dealing with this problem.
The first is to embed a prior in the surrogate design, allowing the model to extrapolate constitutive responses from small datasets.
Examples of these priors are leveraging the geometry or embedding physics into the model acting as an inductive bias~\citep{liu2021review, maia2022physically, klein2024nonlinear, dornheim2024neural}.
These models often have intermediate quantities or additional outputs with a physical meaning that can be leveraged during the training and inference process.

The second approach is to use an active learning strategy instead of only using a pre-trained surrogate.
In active learning, the surrogate is continuously retrained whenever it makes a poor prediction, avoiding the need for an all-encompassing dataset.
Authors in~\cite{rocha2021fly} propose to use the uncertainty of a Gaussian Process (GP) surrogate to inform on-the-fly adaptivity.
There, several fully-solved anchor models are used throughout the simulation to update the GP, triggered by the GP uncertainty.
A similar adaptive approach, $FE^{ANN}$, stores all deformation gradients in a dataset, and by comparing the difference between the invariants of each new deformation to this dataset they can detect when an unseen deformation is encountered~\citep{kalina2023fe}.
When this is the case, they perform new simulations, enrich the dataset, and retrain the neural network.
Adaptively updating the surrogate during the simulation avoids the necessity for all scenarios to be covered in an extensive pre-training dataset.
However, active learning does require the surrogate to be suitable for retraining.
Furthermore, updating the surrogate can cause instabilities or slower convergence of the macroscopic solver.

Instead of using the high-fidelity model just to update the surrogate, a mixture of several models can be applied.
A common strategy is the use of clustering, where elements get clustered based on their response after an initial solve of the system, by classifying the inputs into predetermined clusters~\citep{li2025clustering}. %~\citep{liu2016self}
Alternatively, a mixture of a pre-trained neural network (NN) and a reduced order model (ROM) can be used throughout the domain~\citep{fritzen2019fly}.
After solving an \fetwo{} simulation with the NN surrogate at every point, a second run is performed where all points that reached values outside the NN training region during the first run are replaced with the ROM.
To minimize unnecessary evaluations, the second run can also start with the NN at all points, and switch to the ROM only when strain values outside the NN training region occur.
One of the challenges of these methods is avoiding instability of the finite element method (FEM) caused by stress jumps on boundaries between two models.
One option to avoid these jumps between models is to introduce a transition zone in front of a crack based on the thick level set method~\citep{rocha2020accelerating}.
From a machine learning perspective, these hybrid approaches are analogous to a mixture-of-experts model~\citep{jacobs1991adaptive}, where we use domain-specific information to choose the expert.
Alternatively, this also bears a resemblance to data assimilation methods, such as using a weighted least-squares finite element method to assimilate experimental data with numerical models~\citep{rajaraman2017combining}.

In this work, we allow for adaptively using multiple constitutive models when and where necessary.
Specifically, a probabilistic surrogate and the original model it was trained to replace are aggregated to make predictions.
The main idea is to use the uncertainty in surrogate predictions to guide whether the surrogate or the original model should be used.
For our motivation of accelerating multiscale simulations, the fast surrogate should be used when it has a low uncertainty, and the original model should be used otherwise.
We use the idea of a weighted average transition zone and implement it using a phase field.
The phase field promotes numerical stability, as it naturally introduces a spatial transition zone that is used to switch between the models.
Phase-field models have been used to simulate the evolution of, among others, solidification~\citep{boettinger2002phase}, grain growth~\citep{chen2002phase}, and fracture~\citep{wu2020phase} for several decades now.
In a phase field, a scalar field variable $\phi$ evolves, symbolizing, for example, the transition from a liquid state to a solid state in the context of a two-phase solidification system, or the shift from a solid area to a crack in crack propagation.
Here, we use the phases for the different constitutive models.
The uncertainty coming from the probabilistic surrogate drives the phase field evolution.
At the diffuse interface between phases, we take a weighted average between model predictions to smoothly transition between the models.

Using the phase field brings the cost of solving an additional partial differential equation (PDE) every time step.
In \fetwo{} simulations, the costs associated with solving microstructures are generally so high that the cost of solving the phase field is negligible.
When a surrogate is trained on a well-curated dataset that allows it to predict the full behavior of the simulation with little uncertainty, the phase field remains zero and only the surrogate is used. %, as is currently implicitly assumed by most works using surrogate models,
In this case, this approach only requires the additional computation of the surrogate uncertainty.
The mixture of constitutive models will enable an accurate result in cases where this condition is not met.

The width of the transition zone can play an important role in the stability of the mechanical and phase-field problems.
A narrow transition zone requires fewer expensive model evaluations, making it computationally desirable.
However, this might lead to undesirable stress jumps or make the solver require more Newton-Raphson iterations.
To study this, we compare the phase-field approach to a purely local approach, which does not require solving another PDE.
By letting the uncertainty directly determine which constitutive model is used, we can essentially let the transition zone vanish, and study the impact on the stability of the mechanical problem.

The outline of this paper is as follows.
In Section~\ref{sec:fem_cost}, we discuss using FEM for solving the behavior of structures undergoing loading.
We provide a detailed overview of our phase-field-based approach for mixing models in Section~\ref{sec:hybrid_model}.
There, we also discuss using a GP as a surrogate constitutive model.
We present the results of the different approaches in Section~\ref{sec:results}.
Finally, we conclude with an overview of our findings in Section~\ref{sec:conclusion}.

\section{Finite element analysis}\label{sec:fem_cost}
FEM is used to model the behavior of a structure or material under load.
Specifically, the aim is to find the displacement field $\mathbf{u}$ of a structure subject to boundary conditions.
We obtain equilibrium in the domain by finding the displacement field $\mathbf{u}$ satisfying
\begin{equation}\label{eq:equilibrium}
 \nabla \cdot \boldsymbol{ \sigma } = 0.
\end{equation}
Here, $\boldsymbol{ \sigma }$ is the stress, and $ \nabla \cdot $ indicates the divergence operator.
$\mathbf{u}$ is related to the strains $\boldsymbol{\varepsilon}$ as
\begin{equation}
    \boldsymbol{\varepsilon} = \dfrac{1}{2} \left( \nabla \boldsymbol{u} + \left(\nabla \boldsymbol{u} \right)^T \right),
\end{equation}
and $\boldsymbol{\varepsilon}$ can be related to $\boldsymbol{ \sigma }$ through a constitutive relation:
\begin{equation}\label{eq:const_rel}
    \boldsymbol{\sigma} = \mathcal{C}(\boldsymbol{\varepsilon}, \boldsymbol{\alpha}).
\end{equation}
The material-dependent constitutive operator $\mathcal{C}$ is often nonlinear and can depend on the loading history.
The internal variables $\boldsymbol{\alpha}$ can account for this loading history.
Due to the nonlinear nature of the constitutive model, Equation~\ref{eq:equilibrium} cannot be solved for directly.
Instead, an iterative Newton-Raphson scheme is employed.
To obtain quadratic convergence, we require not just the outcome of the constitutive relation ($\boldsymbol{ \sigma }$), but also its derivative ($\boldsymbol{ D }$).

To find equilibrium, the FEM discretizes the continuous domain of the structure into a number of elements.
The boundary conditions, in the form of forces or prescribed displacements, are applied to the nodes of these elements.
Within each element, the displacement field $\mathbf{u}$ is approximated using simple \textit{basis functions}, often polynomials.
The Gauss quadrature rule is used to evaluate the integrals arising from the discretization, with strategically placed integration points (IPs) within each element.
The constitutive relation from Equation~\ref{eq:const_rel} is evaluated for each integration point.

We consider quasi-static loading, where the load is applied in consecutive time steps without causing dynamic effects.
This allows for obtaining the full load path of the structure of interest.
For a single simulation, we thus have to compute the constitutive model for all quadrature points in the elements, for each Newton-Raphson iteration of all loading steps.

For certain materials, such as composites, the constitutive relation can depend on the geometry at the microscale.
Finding accurate constitutive relations for these materials is challenging, even if the constitutive relations of each constituent are known.
In multiscale simulations using \fetwo{}, a microscopic FEM problem is solved instead, where $\boldsymbol{\varepsilon}$ is imposed as boundary conditions on a microstructure.
We can solve this system as we would in FEM since each of the constituent's properties are well described, and then pass the homogenized stress $\boldsymbol{\sigma}$ back to the macroscale.
However, solving a boundary value problem for every macroscopic quadrature point for every time step brings a considerable computational cost.
This cost is so high that even with parallelization, it sees only limited use in practice, which motivates exploring acceleration strategies for \fetwo{}.

To make it computationally feasible to study the approach in detail, we limit ourselves to single-scale problems in this work.
We use an elasto-plastic material with a von Mises yield criterion as our high-fidelity (HF) model, which serves as the ground truth.
To quantify the acceleration of our approach, it thus makes little sense to compare the computational time --- as evaluating this elasto-plasticity model is about as fast as the surrogate we use to replace it.
Instead, we measure the reduction in the number of HF model evaluations.
The simplification of using an analytical constitutive model therefore does not affect our upcoming conclusions about the approach.

For complex problems, the Newton-Raphson solver can fail to converge for a given load increment.
Adaptive load stepping is then required to obtain a solution.
In this work, we use an adaptive load step strategy that stays constant as long as the solver converges.
If it fails to converge, the load step is reduced by a factor $\gamma$.
This is repeated until the solver converges, or until a minimum load step value is reached.
After reaching the minimum value, increased load steps are attempted, up to a maximum value, after which the simulation is terminated and the problem is considered unsolved.
If it convergences, the load step stays constant for one step before gradually increasing back to the initial increment size, now using a factor $\dfrac{1}{\gamma}$.

\section{Mixture of constitutive models}\label{sec:hybrid_model}
\begin{figure}[bp]
\centering
\begin{minipage}[t]{.47\textwidth}
    \centering
    \includegraphics[width=0.7\textwidth]{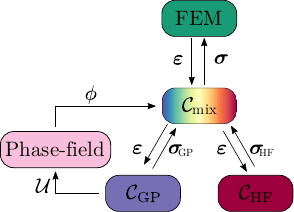}
    \caption{Overview of how the mixture constitutive model interacts with the other models. The tangent stiffness $\boldsymbol{D}$ is omitted for clarity.}
    \label{fig:Mixture_overview}
\end{minipage}
\hfill
\begin{minipage}[t]{.47\textwidth}
    \centering
    \includegraphics[width=.5\textwidth]{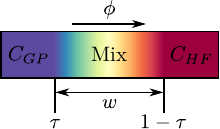}
    \caption{As $\phi$ increases, the constitutive model $\mathcal{C}_{mix}$ gradually switches from $\mathcal{C}_{\mathrm{GP}}$ to $\mathcal{C}_{\mathrm{HF}}$.}
    \label{fig:PF_interface_diagram}
\end{minipage}
\end{figure}

We solve the mechanical problem by combining constitutive models of distinct accuracies and computation speeds in the macroscopic domain.
The analytical elasto-plastic material model serves as the high-fidelity model ($\mathcal{C}_{\mathrm{HF}}$), and we use a GP surrogate model ($\mathcal{C}_{\mathrm{GP}}$).
The general schematic of the approach in Figure~\ref{fig:Mixture_overview} shows what information the mixture model obtains from other models.
Ideally, $\mathcal{C}_{\mathrm{GP}}$ is used throughout most of the domain, and $\mathcal{C}_{\mathrm{HF}}$ only for a small subset of points.
Our hybrid constitutive model is dependent on a phase-field variable $\phi$ to determine which constitutive model is used at each quadrature point:
\begin{equation}\label{eq:hybrid_constitutive}
    \mathcal{C}_{mix} =
    \begin{cases}
    \mathcal{C}_{\mathrm{GP}} & \text{if } \;\;\;\;\;\; \phi < \tau, \\
    (1 - \phi) \mathcal{C}_{\mathrm{GP}} + \phi \mathcal{C}_{\mathrm{HF}} & \text{if } \tau < \phi < 1-\tau, \\
    \mathcal{C}_{\mathrm{HF}} & \text{if } \;\;\;\;\;\; \phi > 1 - \tau.
    \end{cases}
\end{equation}

At the interface between the two phases, the hybrid constitutive model is thus a weighted average of $\mathcal{C}_{\mathrm{HF}}$ and $\mathcal{C}_{\mathrm{GP}}$, as visualized in Figure~\ref{fig:PF_interface_diagram}.
In practice, we take weighted averages over both $\boldsymbol{\sigma}$ and $\boldsymbol{ D }$.
We use a small cutoff value, $\tau$, to avoid simulating $\mathcal{C}_{\mathrm{HF}}$ when the phase field is close to but not exactly zero.
The upper boundary of $1-\tau$ is of lesser importance since the cost of solving $\mathcal{C}_{\mathrm{GP}}$ is negligible, and its uncertainty is still required when updating the mixture.
We use $\tau=0.01$ in all numerical experiments.
Algorithm~\ref{alg:mixture_models} shows this procedure in more detail.
The high-fidelity model generally tracks internal variables $\boldsymbol{\alpha}$ to, for example, account for plasticity.
These internal variables are only committed when moving to the next time step.
When switching from the surrogate model to the full model (when $\phi$ crosses the threshold value $\tau$), the internal variables of the full model are unknown yet necessary.
In this scenario, the history path is recomputed using the full model to obtain the internal variables.
We call this process the re-tracing of the history.
Although this increases the number of evaluations, we only need to re-trace along the converged solutions.
This thus still requires fewer evaluations than if the full model was used for this point from the start.
We could in principle use larger time step sizes when reconstructing, but here we pragmatically opt for retracing every increment missed by the full model.
This also means we have to store converged strain values at every fully-reduced IP for all time steps in the event a retracing is needed.

\begin{algorithm}
\caption{Mixture of constitutive models}\label{alg:mixture_models}
\begin{algorithmic}
\State \textbf{Input: $\boldsymbol{\varepsilon}$, $\phi$}
\State \textbf{Output: $\boldsymbol{\sigma}, \boldsymbol{D}$}
\State \textbf{Initialize:} Internal variables $\boldsymbol{\alpha}$
\If{$\phi < \tau$}
    \State Use only the GP model: $\boldsymbol{\sigma},\boldsymbol{D} = \mathcal{C}_{\mathrm{GP}}(\boldsymbol{\varepsilon})$
\Else
    \State Re-trace $\boldsymbol{\alpha}$ to the current timestep if it is outdated
    \If{$\phi > 1-\tau$}
        \State Use only the high-fidelity model: $\boldsymbol{\sigma},\boldsymbol{D}=\mathcal{C}_{\mathrm{HF}}(\boldsymbol{\varepsilon}, \boldsymbol{\alpha})$
    \Else
        \State $\boldsymbol{\sigma}_{\mathrm{GP}},\boldsymbol{D}_{\mathrm{GP}} = \mathcal{C}_{\mathrm{GP}}(\boldsymbol{\varepsilon})$
        \State $\boldsymbol{\sigma}_{\mathrm{HF}},\boldsymbol{D}_{\mathrm{HF}}=\mathcal{C}_{\mathrm{HF}}(\boldsymbol{\varepsilon}, \boldsymbol{\alpha})$
        \State $\boldsymbol{\sigma} = \phi \boldsymbol{\sigma}_{\mathrm{HF}} + (1-\phi) \boldsymbol{\sigma}_{\mathrm{GP}}  $
        \State $\boldsymbol{D} = \phi \boldsymbol{D}_{\mathrm{HF}} + (1-\phi) \boldsymbol{D}_{\mathrm{GP}}  $
    \EndIf
\EndIf
\end{algorithmic}
\end{algorithm}

\subsection{Phase-field approach}\label{sec:phase-field}
The main purpose of the phase field is to control how the surrogate uncertainty influences the mixture.
Specifically, it allows us to create a controllable interface width and determine the minimum uncertainty required to initiate the switch.
While monolithic phase fields have been used~\citep{gerasimov2016line}, staggered schemes are more common~\citep{luo2023fast} and also what we opt for.
The variational form used in our numerical implementation is as follows:

\begin{equation}\label{eq:pf_variational}
F(\phi, v) =
- \underbrace{\int_\Omega \left[ \mathcal{U} v \right] \, d\Omega}_{\text{driving force}}
+ \underbrace{b \int_\Omega v \, d\Omega}_{\text{opposing force}}
+ \underbrace{\int_\Omega \left[ \epsilon^2 \nabla \phi \cdot \nabla v \right] \, d\Omega}_{\text{interface energy}}
+ \underbrace{\int_\Omega \left[ \omega \phi (1 - \phi)(1 - 2\phi) v \right] \, d\Omega}_{\text{double-well}}
\end{equation}

Here, $v$ represents the test function.
The first term is the driving force $\mathcal{U}$, causing $\phi$ to increase as the uncertainty increases.
Because we solve the phase-field and mechanical problem in a staggered approach, the driving force does not necessarily have to be continuous with respect to the strains to maintain stability.
The second term is a constant opposing force, acting as a threshold that the driving force needs to overcome.
The third term represents the interface energy, penalizing sharp gradients in the phase-field variable.
By influencing the interface width, this term may affect the stability of the problem.
The fourth term is the derivative of the double-well potential, ensuring the phase-field variable favors values near 0 or 1.
The phase-field formulation does not depend on its previous state, as we do not require our phase field to be smooth in time.
This avoids issues related to path dependency, where a different loading step size would lead to different behavior.

\subsection{Staggered updating scheme}\label{sec:NR_schemes}
\begin{figure}[b]
\centering
\includegraphics[width=1.0\textwidth]{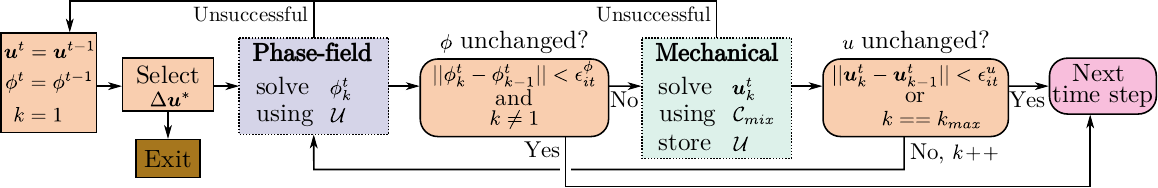}
    \caption{The staggered updating approach used in this work. A prescribed displacement $\Delta u^*$ is selected using an adaptive time-step function, before updating the phase field based on $\mathcal{U}$. Then the mechanical problem is solved. The phase-field and mechanical problem are updated iteratively until convergence, or $k_{max}$ is reached.}
    \label{fig:PF_updating_scheme}
\end{figure}

We solve the phase-field and mechanical problem in a decoupled, staggered manner.
In Figure~\ref{fig:PF_updating_scheme} a schematic overview of the staggered updating scheme is shown.
Although either field can be updated first~\citep{luo2023fast}, we choose to first update the phase field based on the last known surrogate model uncertainty.
The phase field is then fixed while solving the mechanical problem.

Solving each problem only once per time step is computationally desirable, but comes with two problems.
First, the uncertainty from the previous converged step might not be a good indicator of the performance at the current step.
Second, this can introduce a load step size dependency, where taking different step sizes leads to different solutions.
To overcome these issues, we iterate several times over this staggered scheme before moving on to the next step.
We consider this approach to be converged when the norm of the change in $\boldsymbol{u}$ is below a threshold value $\epsilon^u_{it}$.
In certain scenarios, such as for the first few time steps when $\phi=0$ for all points, the phase field will not change, even when $\boldsymbol{u}$ changes.
We then know that $\boldsymbol{u}$ will not change in the next iteration and can consider the current solution converged based on $\phi$, avoiding redundant iterations.
By employing these iterations over the staggered approach, we are adding additional solves of the mechanical problem in which we generally need to perform costly $\mathcal{C}_{\mathrm{HF}}$ evaluations.
However, by using the result of each iteration as the initial prediction for the subsequent iterations, these should converge using fewer Newton-Raphson iterations than they would if starting from $\boldsymbol{u}_{t-1}$.

Instead of relying on the uncertainty coming from the mixed model (which entails expensive high-fidelity simulations), one might instead be inclined to obtain the uncertainty from a preliminary surrogate-only solution.
One would first solve the mechanical problem using only the surrogate model, update the phase field based on the resulting uncertainty, and then solve the mechanical problem again with the mixture model.
However, we find empirically that this approach is not viable: firstly, it heavily depends on reasonable surrogate performance and can fail to converge for inaccurate predictions; secondly, it reintroduces load step size dependency.
Due to these reliability concerns, we do not present results using this approach.

\subsection{Gaussian Processes as surrogate constitutive models}\label{sec:GP}
Bayesian surrogate models provide powerful tools for approximating complex functions.
While several surrogate models provide uncertainty estimates (neural network ensembles, dropout networks~\citep{gal2016dropout}, bootstrapping approaches) and are therefore suitable for our approach, we focus on GP regression.
In~\ref{sec:Gaussianprocess} we provide a general description of GPs for regression that serves as background to this section.
Here we describe how we use GPs specifically as constitutive models.

We create a separate GP for each stress component.
Since our experiments are 2-dimensional, we thus have separate GPs for $\sigma_x$, $\sigma_y$, and $\sigma_{xy}$, each taking the full strain tensor as input.
The GPs use a radial basis function kernel, and their hyperparameters are individually tuned.
Rather than predicting $\sigma$ directly, we instead predict a correction term to a linear elastic model:
\begin{equation}
	\sigma_i = [\boldsymbol{D}_e \boldsymbol{\varepsilon}]_i + \mathrm{GP}_i( \boldsymbol{\varepsilon} ),
\end{equation}
where $\boldsymbol{D}_e$ is the elasticity tensor.
For predictions away from the data, the output thus reverts back to a linear elastic model, rather than to the zero prior.
The scalar driving force for the phase field is computed as the maximum variance of the components:
\begin{equation}
\mathcal{U} = \underset{i}{\max}( \sqrt{var[\mathrm{GP}_i(\boldsymbol{\varepsilon})]} ).
\end{equation}
Since this GP formulation only depends on $\boldsymbol{\varepsilon}$, it is unable to capture elastic unloading.
If unloading is to occur in the simulation, it is likely that the GP could make inaccurate predictions.
For this reason, we focus on scenarios where the global load prescribed by the boundary conditions monotonically increases throughout the simulation.

%%%%%%%%%%%%%%%%%%%%%%%%%%%%%%%%%%%%%%%%%%%%%%%%%%%%%%%%%%%%%%%%%%%%%%%%
\section{Results}\label{sec:results}
We start by showing the behavior of our mixture of models on a simple dogbone structure loaded with a prescribed displacement in tension.
In this dogbone experiment, we focus on the influence of performing several staggered iterations.
Then, moving to a more complex study of a notched plate, we study the influence of the different phase-field parameters.
We additionally compare the phase-field method to an alternative local approach.
Finally, we investigate the potential of our approach to reduce the number of high-fidelity simulations on a plate with holes where a complex phase field is required.

\subsection{Numerical setup}
As discussed earlier, for computational feasibility we do not perform~\fetwo{} in this study, but instead use an analytical constitutive model as $\mathcal{C}_{\mathrm{HF}}$.
For this reason, we do not compare the computational time of our approach to the full~\fetwo{} model, but instead focus on the reduction in the number of high-fidelity model evaluations.
We use an elasto-plastic material with von Mises plasticity in a plane-stress condition, a Young's modulus of $3130$ [MPa], a Poisson's ratio of $0.37$ [-], and yield criterion $\sigma_y = 64.80 - 33.60 \cdot e^{\frac{\varepsilon^p_{eq}}{-0.003407}}$, where $\varepsilon^p_{eq}$ is the equivalent plastic strain (corresponding to the internal variable $\boldsymbol{\alpha}$ used earlier).
These material properties are adopted from~\citep{rocha2019numerical}.
In all experiments, we use quadratic triangular elements (T6) for the mechanical problem and linear triangular elements (T3) for the phase field.
For the dogbone structure, a constant displacement increment of $0.001$ is sufficient to obtain convergence.
In the remaining problems, the adaptive time step handling as described in Section~\ref{sec:fem_cost} is used, with $\gamma=0.5$.

We study the influence of the phase field with three different GP surrogates, \textit{GP10}, \textit{GP30}, and \textit{GP100}, trained on different datasets with 10, 30, and 100 load curves, respectively.
Each load curve consists of 20 monotonically increasing steps in a random strain direction, up to a strain norm of 10\%.
We choose 10\% for pragmatic reasons, assuming we have no prior information about the strains that will occur during the simulation.
Just because certain GPs have more data does not guarantee that their predictions during the simulation are always more accurate, since the loading directions of their data are random.
Still, it can be reasonably expected that the GPs with more data improve accuracy and reduce reliance on the full model.

The framework is implemented using the \textit{FEniCSx} finite element library~\citep{dolfinx2023} with the \textit{dolfinx\_materials}~\citep{bleyer2024dolfinx_materials} package to define arbitrary constitutive models.
\textit{PETSc}~\citep{petsc-web-page} solvers are used to solve the mechanical and phase-field problem.
The phase field is solved using a constrained solver with bounds $[0,1]$.
All code is available on GitHub at \url{https://github.com/SLIMM-Lab/phase-field-mixture}.

\subsection{Time step consistency}\label{sec:res_updating_schemes}
\begin{figure}[!htbp]
    \centering
    \begin{minipage}[b]{0.49\textwidth}
        \centering
        \includegraphics[width=\textwidth]{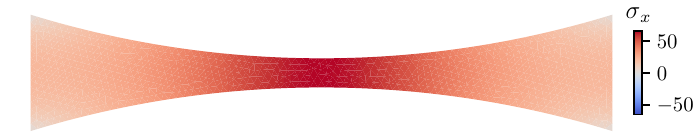}
        \subcaption{$\sigma_{xx}$}
    \end{minipage}
    \begin{minipage}[b]{0.49\textwidth}
        \centering
        \includegraphics[width=\textwidth]{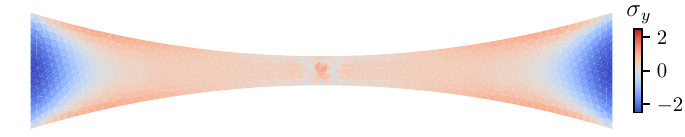}
        \subcaption{$\sigma_{yy}$}
    \end{minipage}

    \vspace{0.1cm}

    \begin{minipage}[b]{0.49\textwidth}
        \centering
        \includegraphics[width=\textwidth]{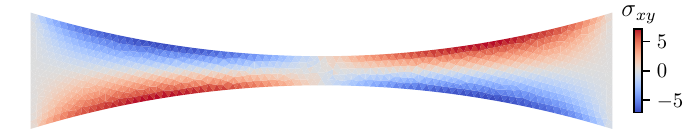}
        \subcaption{$\sigma_{xy}$}
    \end{minipage}
    \begin{minipage}[b]{0.49\textwidth}
        \centering
        \includegraphics[width=\textwidth]{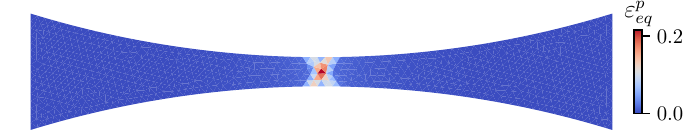}
        \subcaption{$\varepsilon^p_{eq}$}
    \end{minipage}

    \caption{Full-field stress and equivalent plastic strain $\varepsilon^p_{eq}$ plots of the dogbone experiment at the maximum prescribed displacement.}
    \label{fig:dogbone_true}
\end{figure}

We start with a simple dogbone structure to study the influence of the staggered approach of Figure~\ref{fig:PF_updating_scheme} on the time-step dependency.
The dogbone structure is fixed on its left edge, and loaded by incrementally applying a prescribed displacement to the right edge, up to an elongation of 2\%.
Figure~\ref{fig:dogbone_true} shows the true full-field responses for the stresses and equivalent plastic strain.
The stresses mainly occur in the x-direction and concentrate at the center of the dogbone, where the structure is narrower than at its ends, with ensuing strain localization.

\begin{figure}
\centering
\includegraphics[width=.92\textwidth]{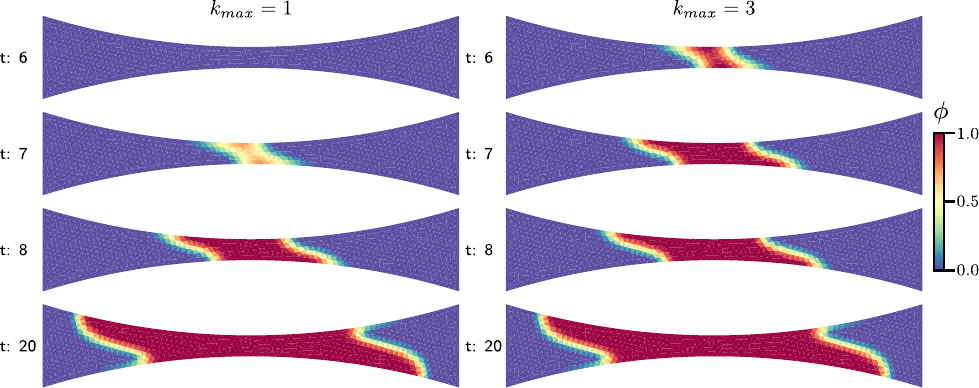}
    \caption{Evolution of the phase field for the dogbone experiment. The phase field gradually grows until time step 17, after which it stays constant due to strain localization. The asymmetrical pattern occurs as a result of the asymmetrical shear stresses.}
    \label{fig:PF_evolution_dogbone_fullfield}
\end{figure}

For the phase-field mixture we use the \textit{GP10} surrogate with hyperparameters $\epsilon=10^{-2}$, $\omega=10^{-3}$, and $b=1$.
On the left side of Figure~\ref{fig:PF_evolution_dogbone_fullfield}, we plot the phase-field variable for various time steps using $k_{max}=1$ (i.e. a single solve of each PDE per load increment, see Figure~\ref{fig:PF_updating_scheme}).
At the start, when $\phi$ is zero everywhere, only the GP surrogate is used.
Starting from time step seven, the uncertainty in the high-strain region at the center of the dogbone drives the phase field locally up to one.
The phase field can switch from $\phi=0$ to $\phi=1$ in a single load step when the uncertainty increases quickly.
This is necessary to prevent large errors when the GP suddenly becomes uncertain.

\begin{figure}[htb]
\centering
\includegraphics[width=0.4\textwidth]{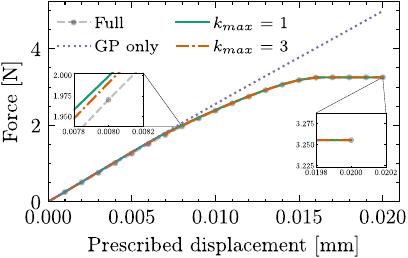}
    \caption{The load-displacement curves for the dogbone experiment. While using only the surrogate model results in poor performance, the hybrid approach shows a large agreement with the baseline full model, with a slightly lower error for $k_{max}=3$ compared to $k_{max}=1$. }
    \label{fig:dogbone_FU_comparison}
\end{figure}

In Figure~\ref{fig:dogbone_FU_comparison}, we present the load-displacement curves for this problem.
Using only the GP surrogate results in a poor prediction that does not capture the plastic behavior of the material, demonstrating that the surrogate model is inadequate to predict the structure by itself.
Still, by switching away from the surrogate at the right time, the hybrid approach closely follows the full model path.

\begin{figure}[htbp]
\centering
\includegraphics[width=.4\textwidth]{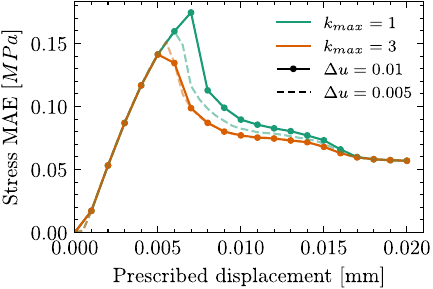}
    \caption{The full-field stress error when using $k_{max}=1$ or $k_{max}=3$ for different load step increments.}
    \label{fig:Timestep_consistency}
\end{figure}

Figure~\ref{fig:Timestep_consistency} shows the mean absolute error of this hybrid model over the domain compared to using the high-fidelity model everywhere.
The errors are computed for two fixed increment sizes $\Delta u$.
For $k_{max}=1$, we observe that when the phase field is changing, from $u=0.07$ to $u=0.17$, there is a clear difference in errors when using different increments.
This difference can be explained by the phase field using uncertainty from the previous time step to determine the mixture of constitutive models.
When using a different increment size, the phase field thus updates on different information, creating a time-step dependency.
In contrast, we observe that using $k_{max}=3$ results in the same response for both loading increment sizes at the points where both are evaluated.
Performing multiple iterations of the staggered approach per loading increment thus avoids the time-step dependency.

\begin{figure}
\centering
\begin{minipage}{.47\textwidth}
    \centering
\includegraphics[width=.8\textwidth]{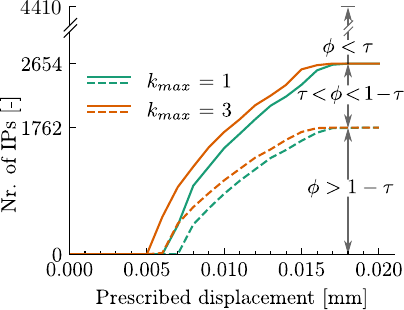}
    \caption{Number of integration points (IPs) that use the GP surrogate ($\phi < \tau$), a mixture of both models ($\tau < \phi < 1-\tau$), or the high-fidelity model only ($\phi>1-\tau$). The mesh contains 4410 IPs in total.}
    \label{fig:PF_dogbone_num_elements}
\end{minipage}
\hfill
\begin{minipage}{.47\textwidth}
    \centering
\includegraphics[width=.8\textwidth]{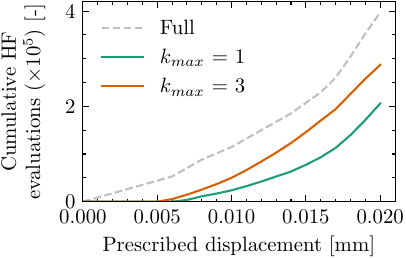}
    \caption{Comparison of the cumulative number of HF evaluations. Retraced history computations are counted in the time step in which they are retraced. The significant difference between $k_{max}=1$ and $k_{max}=3$ comes mainly from the additional staggered iterations, requiring more Newton-Raphson steps.}
    \label{fig:PF_dogbone_num_elements_cumulative}
\end{minipage}
\end{figure}

By performing multiple iterations, the phase field is updated based on information from the current load increment, rather than based on that of the last iteration.
In Figure~\ref{fig:PF_dogbone_num_elements} the evolution of $\phi$ throughout the simulation is plotted.
It is clear from this plot that $k_{max}=1$ lags roughly one step behind compared to $k_{max}=3$.
This can also be observed by comparing the left and right sides of Figure~\ref{fig:PF_evolution_dogbone_fullfield}.

The downside of performing several iterations of the staggered approach is that it requires more Newton-Raphson iterations of the mechanical problem, causing more high-fidelity constitutive model evaluations.
Figure~\ref{fig:PF_dogbone_num_elements_cumulative} shows the cumulative number of high-fidelity evaluations, where we observe that $k_{max}=3$ uses more evaluations than $k_{max}=1$.
The reduction in high-fidelity evaluations compared to running the full model is relatively small in this example, as we intentionally use a GP trained on a very limited dataset to show the behavior of the phase field.
Since it avoids the time-step-dependent behavior and improves the model's accuracy, we argue that performing several iterations is worthwhile for guaranteeing model robustness.

\begin{figure}
\centering
\includegraphics[width=.4\textwidth]{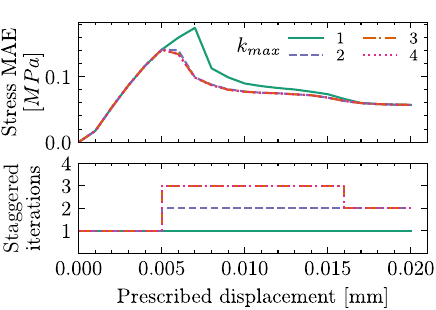}
    \caption{The influence of $k_{max}$ on the simulation. In this example, the model always converges within three iterations.}
    \label{fig:consistent_comparison}
\end{figure}

We further investigate the influence of the number of staggered iterations by running the model for more values of $k_{max}$.
Figure~\ref{fig:consistent_comparison} shows the results by plotting the error and number of staggered iterations used throughout the simulation.
Using two iterations already shows a significant difference in the error compared to using only one.
Even with $k_{max}=4$, the model always converges within three iterations for this experiment.
However, it is possible for the problem to oscillate between iterations, never meeting the tolerances, so $k_{max}$ should be set.
We set $k_{max}=3$ for the remainder of this work.

\subsection{Phase-field parameters}
In this section, we explore the influence of the phase-field parameters in more detail.
To highlight their influence, we create an experimental setup with a more complex phase-field evolution pattern.

\subsubsection{Experimental setup}
\begin{figure}[b]
\centering
\includegraphics[width=.6\textwidth]{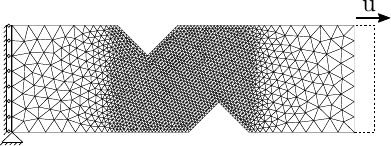}
    \caption{Overview of the notched plate example.}
    \label{fig:notched_plate}
\end{figure}

\begin{figure}[t]
    \centering
    \begin{minipage}[b]{0.49\textwidth}
        \centering
        \includegraphics[width=\textwidth]{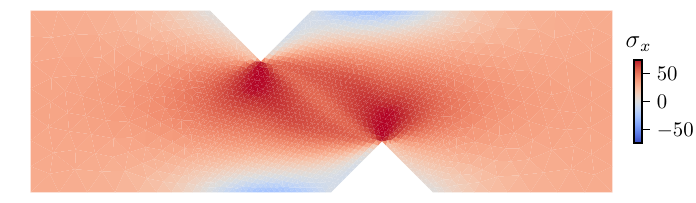}
        \subcaption{$\sigma_{xx}$}
    \end{minipage}
    \begin{minipage}[b]{0.49\textwidth}
        \centering
        \includegraphics[width=\textwidth]{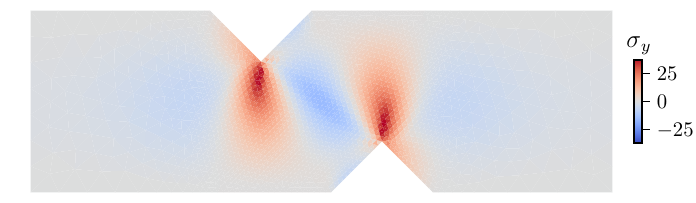}
        \subcaption{$\sigma_{yy}$}
    \end{minipage}

    \vspace{0.1cm}

    \begin{minipage}[b]{0.49\textwidth}
        \centering
        \includegraphics[width=\textwidth]{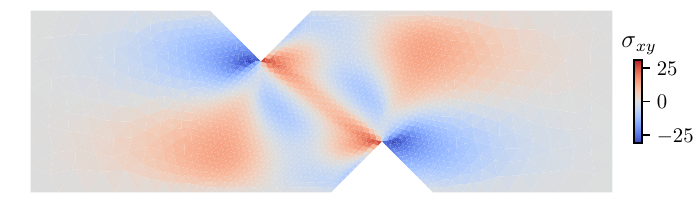}
        \subcaption{$\sigma_{xy}$}
    \end{minipage}
    \begin{minipage}[b]{0.49\textwidth}
        \centering
        \includegraphics[width=\textwidth]{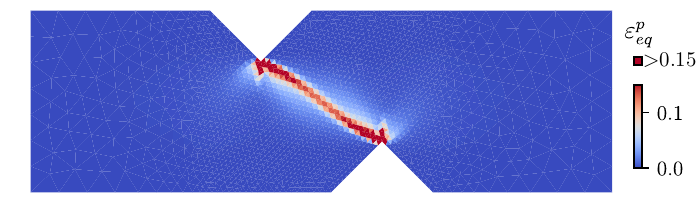}
        \subcaption{$\varepsilon^p_{eq}$}
    \end{minipage}

    \caption{Full-field stress and equivalent plastic strain $\varepsilon^p_{eq}$ plots of the notched plate experiment at the end of the prescribed displacement.}
    \label{fig:notched_plate_true}
\end{figure}

A plate with two notches diagonally opposite each other is loaded by incrementally displacing the right boundary by a prescribed amount.
This setup is visualized in Figure~\ref{fig:notched_plate}, and we plot full-order simulation results in Figure~\ref{fig:notched_plate_true}.
It can be observed that stress concentrations occur at the tips of the notches, and a plastic strain band develops between these tips.
We vary the element size over the domain to capture the strain localization accurately.
Even though the global loading increases monotonically, we observe local unloading in some integration points close to the notch tips.
As our GP surrogates cannot capture unloading, switching to the high-fidelity model in these points becomes essential.
This is, however, not guaranteed by our method, particularly when using a high opposing force (through parameter $b$) that allows the GP to be used despite high uncertainty.
Still, we empirically find that the limited local unloading has minimal impact on the results and does not affect our conclusions.

Several quantities are used to compare the influence of the hyperparameters.
To evaluate the accuracy of the model, we compare force-displacement curves.
Since we now use an adaptive stepper, the points of the force-displacement curve do not coincide if a certain step fails to converge, so we cannot directly compute their differences.
Instead, this error is computed by first linearly interpolating each F-u curve to the fixed displacement points that the solution follows if it always converges.
To compute the error, we sum the absolute force differences between the approximate model and the true model.
In addition, we plot the average number of integration points (IPs) where the high-fidelity model is used (i.e. where $\phi > \tau$).
This indicates how much the high-fidelity phase has grown during the simulation.
As an indication of the stability, we compare the required number of Newton-Raphson iterations of the mechanical problem, obtained by summing over all adaptive time steps (converged or unconverged), for all staggered approach loops.
For the number of high-fidelity evaluations, we track each time the high-fidelity model is evaluated, including the re-tracing of load paths to obtain the internal variables (as described in Section~\ref{sec:NR_schemes}).

\subsubsection{Phase-field opposing force}\label{sec:pf_driving_force}
We study the influence of the phase-field hyperparameter $b$ on the mechanical problem.
The parameter $b$ determines the magnitude of the constant opposing force, as shown in Equation~\ref{eq:pf_variational}.
The larger $b$ is, the larger the uncertainty $\mathcal{U}$ needs to be to drive the phase field up from zero towards one.
Since the uncertainty of different surrogate models can behave very differently, it is challenging to directly compare them for the same $b$ value.
For example, the GPs we use in this work all have different hyperparameters that influence their maximum variance and how quickly the uncertainty increases away from training points.
While this maximum variance could in principle be used to scale $\mathcal{U}$, making $b$ an adimensional threshold, other surrogate methods often do not have an upper bound to their uncertainty.
Therefore, we let $b$ simply be the unscaled uncertainty.
\begin{figure}[htbp]
\centering
\includegraphics[width=.8\textwidth]{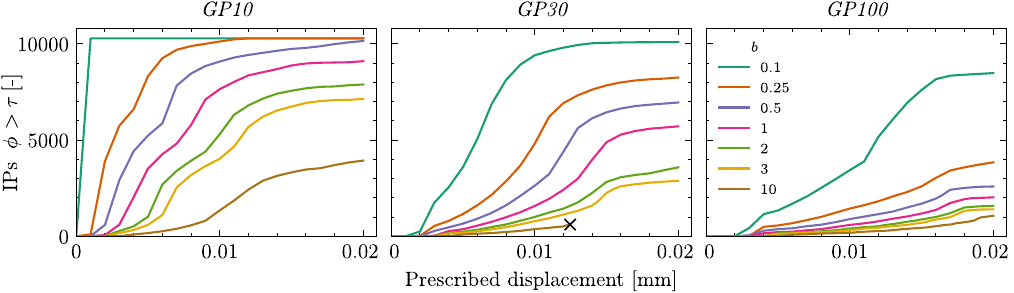}
    \caption{Evolution of the phase field during the simulation for various $b$ values. Default values $\epsilon=10^{-2}$ and $\omega=10^{-3}$ are used.
    The simulation with \textit{GP30}, $b=10$, failed to converge.
    }
    \label{fig:b_influence_phi_evolution}
\end{figure}

\begin{figure}[tbp]
\centering
\includegraphics[width=.32\textwidth]{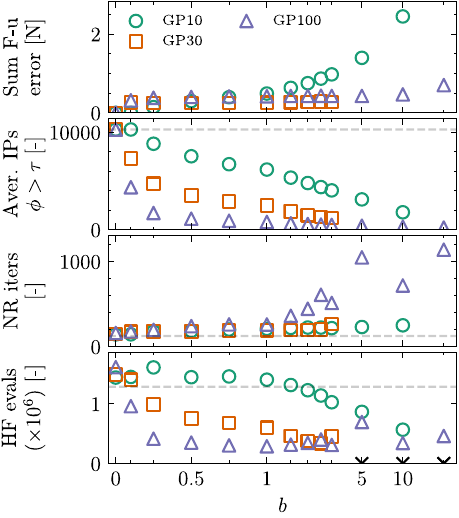}
    \caption{Influence of the opposing force $b$ on the mechanical problem. Note that the x-axis is linear from zero to one, and then scales logarithmically.}
    \label{fig:b_influence}
\end{figure}

We vary $b$ while fixing the other phase field parameters from Equation~\ref{eq:pf_variational} to $\epsilon=10^{-2}$ and $\omega=10^{-3}$.
In Figure~\ref{fig:b_influence_phi_evolution} we show how the phase field evolves throughout the simulation for the different GPs for different values of $b$.
We observe that the uncertainty of GPs with more data grows more slowly and that the larger the value of $b$, the fewer IPs use the high-fidelity model.
We visualize the influence of $b$ on the error, the phase-field evolution, the required number of Newton-Raphson iterations, and the number of high-fidelity evaluations in Figure~\ref{fig:b_influence}.
The gray line in the plots refers to the values obtained when running the full model, with the high-fidelity constitutive model used in all IPs.
Missing data points, and the crosses in the bottom plot, indicate that the setting failed to converge (black crosses indicate a failure of the mechanical problem, red crosses indicate a failure of the phase field).
The error of the \textit{GP10} model increases consistently as we increase $b$.
This indicates that we are using the insufficiently trained GP model even when it is highly uncertain.
As expected, the GP models with more training data result in fewer $\phi > \tau$ elements.
As we increase $b$, the required number of Newton-Raphson iterations of the mechanical problem generally increases, with several instances of \textit{GP10} and \textit{GP30} failing to reach the end of the simulation.
This creates a trade-off where a higher $b$ requires less high-fidelity evaluations but can increase the error and lead to instabilities.

\subsubsection{Phase-field double-well}\label{sec:pf_omega}
\begin{figure}[tbp]
\centering
\includegraphics[width=.32\textwidth]{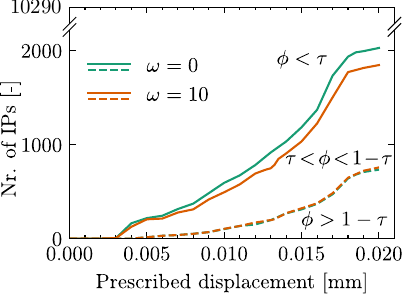}
    \caption{Influence of the double-well term $\omega$ on the phase-field values. A higher $\omega$ leads to more elements having a $\phi$ value of either below $\tau$ or above $1-\tau$. We use \textit{GP100}, $\epsilon=10^{-2}$, and $b=1$.}
    \label{fig:omega_influence_elements}
\end{figure}

In Figure~\ref{fig:omega_influence_elements} we show the evolution of the phase-field values throughout a simulation for $\omega=0$ and $\omega=10$.
We observe that this double-well term of the phase field pushes values away from $\phi=0.5$ and towards $\phi=0$ and $\phi=1$.
Increasing $\omega$ thus has the potential to increase the number of IPs with $\phi < \tau$, avoiding the high-fidelity model evaluation.
Having fewer elements at the interface reduces the number of high-fidelity evaluations, leading to a larger acceleration.
However, by setting the parameter $\omega$ too high, solving the phase field itself can become challenging.
The difference for the $\phi = 1-\tau$ boundary is smaller, indicating that there are fewer points with $\phi$ values close to that boundary.

\begin{figure}[tbp]
\centering
\includegraphics[width=.32\textwidth]{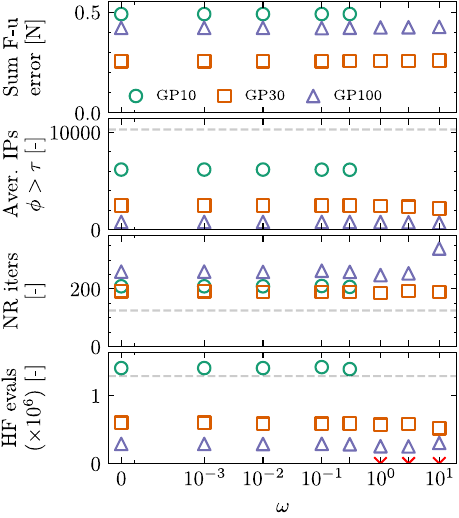}
    \caption{Influence of the double-well term $\omega$ on the mechanical problem.}
    \label{fig:omega_influence}
\end{figure}

The influence of $\omega$ on the result of the simulation for different GPs is shown in Figure~\ref{fig:omega_influence}.
For \textit{GP10}, the phase field reaches the part of the domain with larger mesh elements.
The combination of \textit{GP10} being inaccurate and $\omega$ being high makes convergence more difficult.
Values below $\omega<10^{-1}$ appear to have little effect on any of the quantities.
For larger values of $\omega$, slightly fewer $\mathcal{C}_{\mathrm{HF}}$ evaluations are required, yet the phase field also becomes unstable, failing to converge for some cases, or requiring significantly more Newton-Raphson iterations.
As the impact of this double-well term appears small, we conservatively set $\omega=10^{-3}$ for the remainder of this work.

\subsubsection{Phase-field interface width}\label{sec:pf_epsilon}
\begin{figure}
\centering
\includegraphics[width=.32\textwidth]{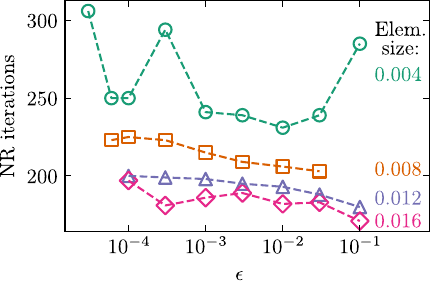}
    \caption{Stability of the phase-field and mechanical problem concerning the element size and phase-field variable $\epsilon$ on the notched plate problem. All problems are evaluated on an $\epsilon \in [10^{-5}, 10^{-1}]$ interval, but only converged points are shown. We have used $\omega=10^{-3}$ and the \textit{GP10} model for all simulations.}
    \label{fig:epsilon_mesh_size}
\end{figure}

\begin{figure}[!htbp]
    \centering
    \hfill
    \begin{minipage}[b]{0.32\textwidth}
        \centering
        \includegraphics[width=\textwidth]{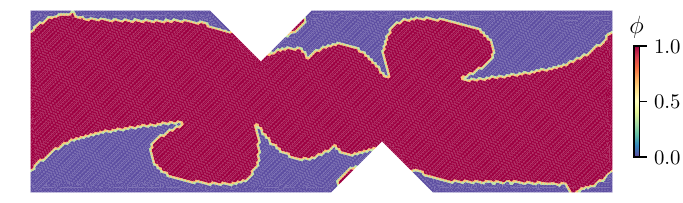}
        \subcaption{Element size: 0.004. $\epsilon=0.00003$}
    \end{minipage}
    \begin{minipage}[b]{0.32\textwidth}
        \centering
        \includegraphics[width=\textwidth]{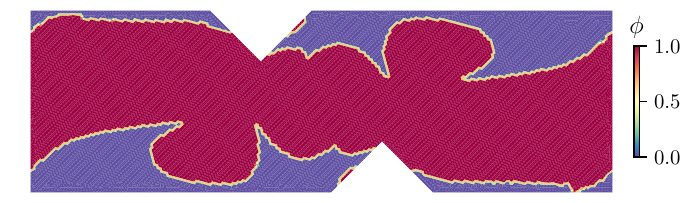}
        \subcaption{Element size: 0.004. $\epsilon=0.0001$}
    \end{minipage}
    \begin{minipage}[b]{0.32\textwidth}
        \centering
        \includegraphics[width=\textwidth]{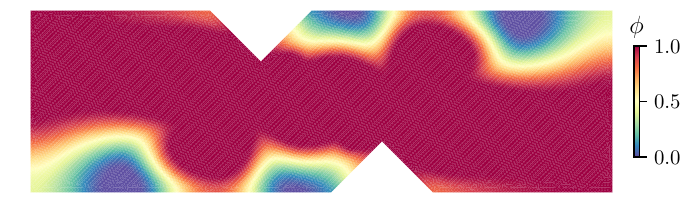}
        \subcaption{Element size: 0.004. $\epsilon=0.03$}
    \end{minipage}

    \vspace{0.3cm}

    \hfill
    \begin{minipage}[b]{0.32\textwidth}
        \centering
        \includegraphics[width=\textwidth]{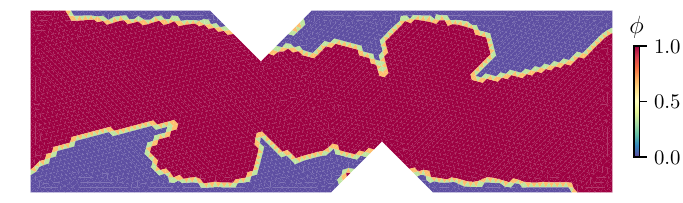}
        \subcaption{Element size: 0.008. $\epsilon=0.0001$}
    \end{minipage}
    \begin{minipage}[b]{0.32\textwidth}
        \centering
        \includegraphics[width=\textwidth]{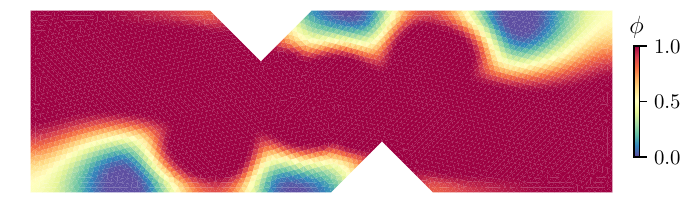}
        \subcaption{Element size: 0.008. $\epsilon=0.03$}
    \end{minipage}

    \vspace{0.3cm}

    \hfill
    \begin{minipage}[b]{0.32\textwidth}
        \centering
        \includegraphics[width=\textwidth]{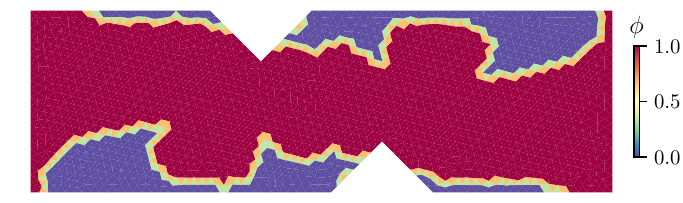}
        \subcaption{Element size: 0.012. $\epsilon=0.0001$}
    \end{minipage}
    \begin{minipage}[b]{0.32\textwidth}
        \centering
        \includegraphics[width=\textwidth]{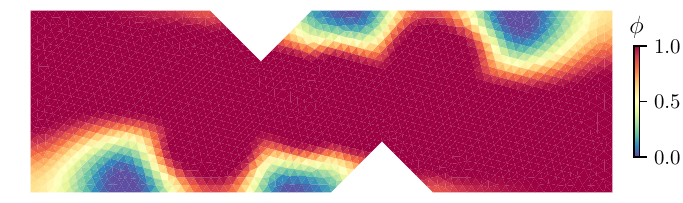}
        \subcaption{Element size: 0.012. $\epsilon=0.03$}
    \end{minipage}

        \vspace{0.3cm}

    \hfill
    \begin{minipage}[b]{0.32\textwidth}
        \centering
        \includegraphics[width=\textwidth]{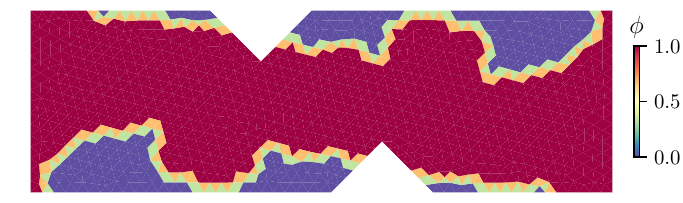}
        \subcaption{Element size: 0.016. $\epsilon=0.0001$}
    \end{minipage}
    \begin{minipage}[b]{0.32\textwidth}
        \centering
        \includegraphics[width=\textwidth]{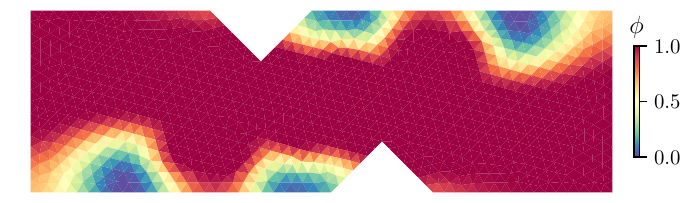}
        \subcaption{Element size: 0.016. $\epsilon=0.03$}
    \end{minipage}

    \caption{The phase-field variable at the final time step of the notched plate problem for various element sizes and values of the phase-field interface parameter $\epsilon$.}
    \label{fig:full_field_epsilon}
\end{figure}

The phase-field parameter $\epsilon$ influences the width of the interface between the phases and, therefore, the width of the transition region between the constitutive models (the interface width is equal to approximately $4 \epsilon$). % and is important for the convergence of the phase-field problem and the mechanical problem.
To study the influence of the interface more carefully, we perform tests with a constant element size throughout the domain.
In Figure~\ref{fig:epsilon_mesh_size}, the influence of $\epsilon$ on the required number of Newton-Raphson iterations is shown for different element sizes.
A few corresponding full-field plots of the phase-field value for the different mesh sizes are plotted in Figure~\ref{fig:full_field_epsilon}.
We can observe two main trends in these figures.
Firstly, finding solutions for small $\epsilon$ values requires smaller mesh sizes.
A minimal interface width is thus required to solve the phase field (but the wider the interface, the more high-fidelity models are evaluated in each iteration).
As we decrease $\epsilon$ such that the mesh element size is smaller than $4 \epsilon$, the interface width stays one element until the phase field fails to converge.
Secondly, the number of Newton-Raphson iterations of the mechanical problem decreases slightly as we increase $\epsilon$.
When the phase field fails for some steps, but the overall simulation still converges with adaptive steps, it causes the spikes observed for element size 0.004.
The general trend indicates that having a transition region span over multiple elements improves the numerical stability of the mechanical problem.

The mesh used for the other notched plate experiments in this work has a variable mesh element size to capture strain localization with fine elements (0.005) while having coarse elements (up to 0.04) in the rest of the domain, as shown in Figure~\ref{fig:notched_plate} and~\ref{fig:notched_plate_true}.
When using such a variable mesh, a separate mesh with a different characteristic element size could be implemented for the phase field as in~\citep{goswami2019adaptive}.
This brings additional challenges and is not considered here.

\begin{figure}
\centering
\includegraphics[width=.32\textwidth]{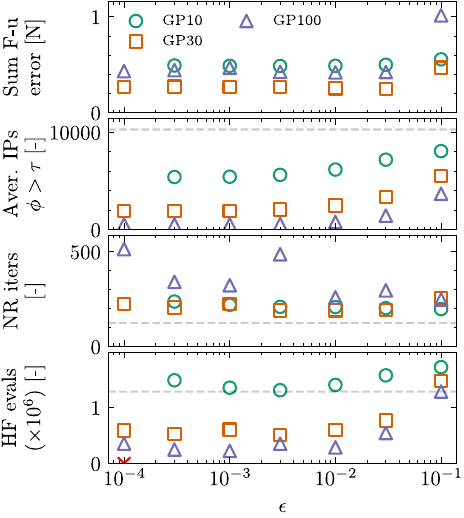}
    \caption{The influence of the phase-field interface parameter $\epsilon$ on the mechanical and phase-field problem.}
    \label{fig:epsilon_influence}
\end{figure}

Figure~\ref{fig:epsilon_influence} shows the influence of $\epsilon$ for the various GP surrogates on the domain with a variable mesh size.
As $\epsilon$ decreases, the interface becomes narrower, leading to fewer elements where $\phi > \tau$.
For $\epsilon=10^{-1}$ the error increases, especially for \textit{GP100}.
This is because when there is a region with only a few points with high uncertainty, a large $\epsilon$ can cause the phase field to not fully switch, hindering accuracy.
The difference in accuracy is generally small, and the number of high-fidelity evaluations reduces for smaller $\epsilon$ values.
Still, as observed earlier, when $\epsilon$ is small, the mechanical and phase-field problems become less stable and require more NR iterations or do not converge.
Therefore, the default value we have used so far of $\epsilon=10^{-2}$ appears to be a sensible choice.

%%%%%%%%%%%%%%%%%%%%%%%%%%%%%%%%%%%%%%%%%%%%%%%%%%%%%%%%%%%%%%%%%%%%%%%%
\subsection{Local mixture}\label{sec:local_mixture_results}
To avoid being limited by the convergence of the phase field for small $\epsilon$ values, we compare the phase field to a simple approach where no phase-field model is solved.
In the previous demonstrations, the phase-field spatial smoothening, a non-local phenomenon, was necessary for the phase field to converge.
To better study the influence of the transition zone without being limited by having to solve the phase-field PDE, we consider an alternative case with no spatial smoothening.
Here, we directly use the driving force to evolve the mixture.
We refer to this as the local approach since the mixture of models in each IP is independent of the uncertainty in other IPs.
In addition to removing the transition zone, another benefit is that we no longer need to solve the phase field, reducing the complexity of the method.
We still use a similar formulation to the phase field, with $b$ opposing the uncertainty.
We consider two variations of the local approach.
The first version linearly shifts between the models:
\begin{equation}
    \phi = \mathcal{U} - b, \quad \phi \in [0,1].
\end{equation}
In contrast to the phase field, this local shift does not necessarily lead to a spatial transition zone.
When two points close in space have very different $\mathcal{U}$ values, the spatial transition can still be sharp.
However, the linear shift ensures that a small change in strain does not cause a completely different constitutive model to be used --- the stress remains continuous with respect to the strain.
Alternatively, we can use a step function, directly switching between the models when the uncertainty reaches the threshold:
\begin{equation}
    \phi =
    \begin{cases}
    0 & \text{if } \mathcal{U} < b, \\
    1 & \text{if } \mathcal{U} \geq b.
    \end{cases}
\end{equation}
Note that both versions require the high-fidelity model when $\mathcal{U} \geq b$.
Similar to the phase-field approach, this method suffers from time-step dependency because the mixture uses the uncertainty of the previous step.
Therefore, we also perform several iterations of the staggered approach, with $k_{max}=3$.

\begin{figure}[htbp]
\centering
\begin{minipage}[t]{.32\textwidth}
    \centering
\includegraphics[width=1.0\textwidth]{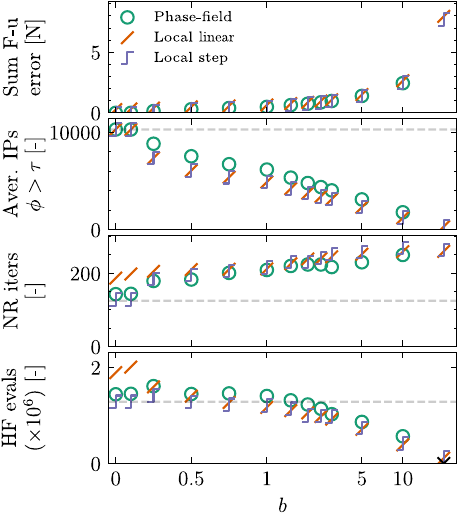}
    \subcaption{GP 10}
    \label{fig:b_GP_comparison_gp10}
\end{minipage}
\hfill
\begin{minipage}[t]{.32\textwidth}
    \centering
\includegraphics[width=1.0\textwidth]{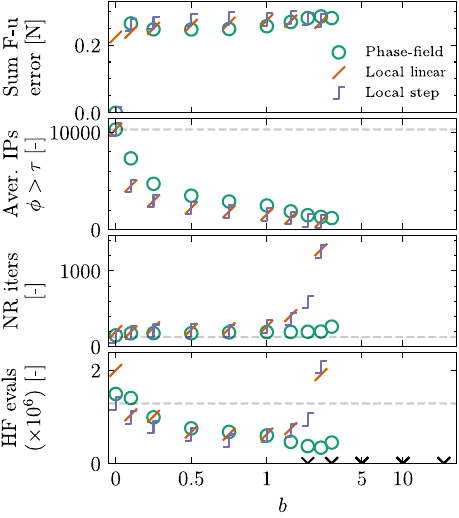}
    \subcaption{GP 30}
    \label{fig:b_GP_comparison_gp30}
\end{minipage}
    \hfill
\begin{minipage}[t]{.32\textwidth}
    \centering
\includegraphics[width=1.0\textwidth]{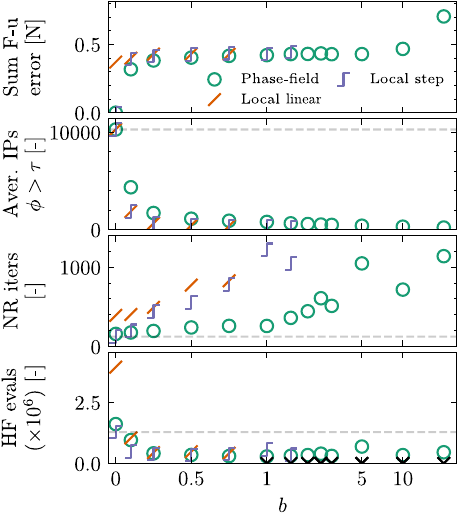}
    \subcaption{GP 100}
    \label{fig:b_GP_comparison_gp100}
\end{minipage}
\caption{Comparison of the opposing force parameter b for the various approaches. For the phase field we set $\epsilon=10^{-2}$ and $\omega=10^{-3}$. Note that the x-axis is linear from zero to one, and then scales logarithmically. All simulations that fail to converge fail due to the mechanical problem.}
\label{fig:b_GP_comparison_total}
\end{figure}

In Figure~\ref{fig:b_GP_comparison_total}, we show how $b$ influences the results for these two local approaches compared to the phase-field approach.
Across all simulations, the errors are nearly identical between the models.
Because the local methods do not have an interface, they always have fewer IPs that require high-fidelity model evaluations for a given uncertainty.
The number of Newton-Raphson iterations is similar when $b$ is low.
However, for larger $b$ values, the local approaches require significantly more iterations with \textit{GP30} and \textit{GP100}.
In addition, the local approaches start failing for these models as $b$ increases.
Removing the interface can thus be beneficial in reducing the number of high-fidelity evaluations for specific settings, yet generally makes the model less stable and less robust.

\subsection{Surrogate data}\label{sec:GP_data}
So far, we have shown the influence of the various parameters on the behavior of the mixture of models.
Now, we study the impact of the training dataset size on the reduction of high-fidelity samples in more detail.
To investigate whether this framework can handle complex phase fields, we consider a plate with holes to create numerous local regions with high stresses that challenge the GP surrogates and induce a switch.
We constrain the left side of the domain and incrementally prescribe a displacement to the right side of the domain.
We use $\epsilon=10^{-2}$, $\omega=10^{-3}$, and $b=1$ for the phase field, we do not consider the local approach here.

\begin{figure}[htbp]
\centering
\includegraphics[width=.32\textwidth]{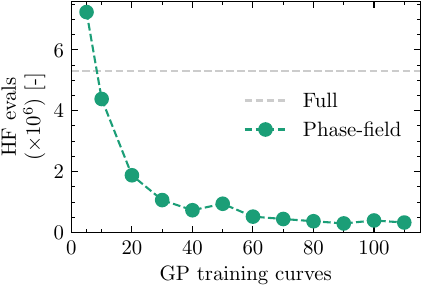}
    \caption{GP training dataset size influence on the number of high-fidelity evaluations. The plot shows the average of 5 runs with different datasets. The full-model is used as a reference. For very small datasets, our approach can use more evaluations than the full model due to the staggered iterations ($k_{max}=3$). }
    \label{fig:gp_dataset_comp}
\end{figure}

\begin{figure}[!htbp]
    \centering
    \begin{minipage}[b]{0.32\textwidth}
        \centering
        \includegraphics[width=\textwidth]{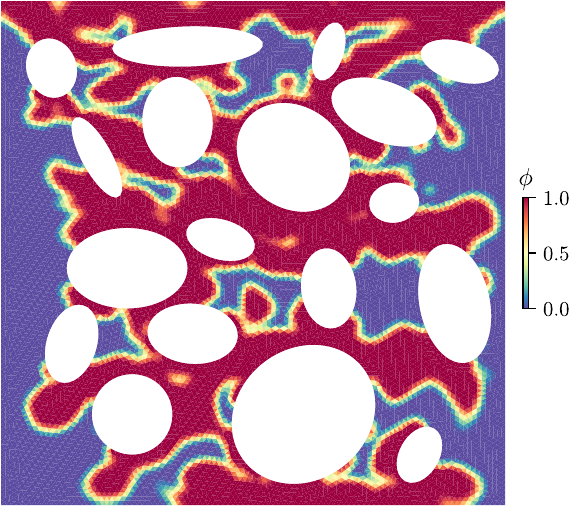}
        \subcaption{\textit{GP10}}
    \end{minipage}
    \begin{minipage}[b]{0.32\textwidth}
        \centering
        \includegraphics[width=\textwidth]{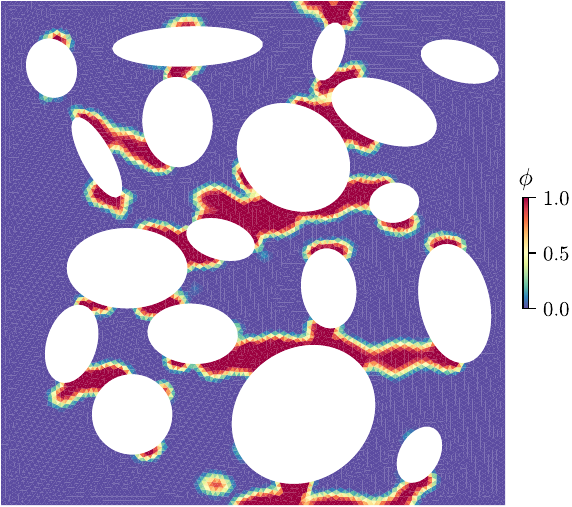}
        \subcaption{\textit{GP30}}
    \end{minipage}
    \begin{minipage}[b]{0.32\textwidth}
        \centering
        \includegraphics[width=\textwidth]{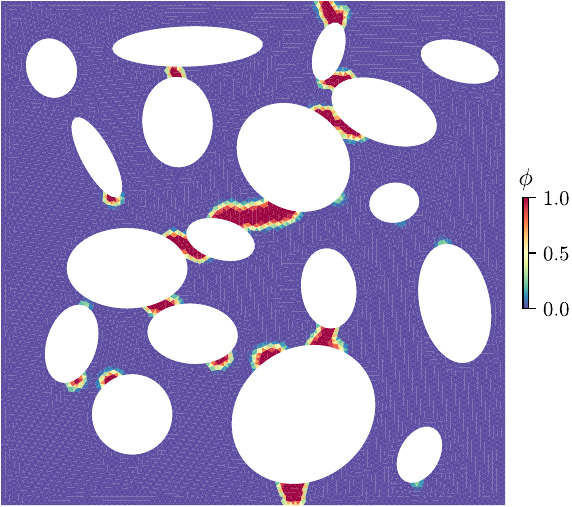}
        \subcaption{\textit{GP100}}
    \end{minipage}

    \caption{Phase-field variable $\phi$ at the maximum prescribed displacement for surrogates with varying amounts of training data.}
    \label{fig:fullfield_plateholes}
\end{figure}

Figure~\ref{fig:gp_dataset_comp} shows the number of high-fidelity evaluations during the simulation as we vary the amount of training data in the GP surrogate.
For three of these cases, we plot the full-field phase field at the end of the simulation in Figure~\ref{fig:fullfield_plateholes}.
These figures show the number of high-fidelity evaluations quickly decreasing as the GP is trained with more training data.
Still, at a certain point, the downward trend appears to converge, and the number of high-fidelity evaluations is not further reduced.
This is because the GP fails to capture localization with strains outside its training data range, which go up to a strain norm of 10\%, a choice that must be made prior to running the simulation.
This shows the difficulty of creating an appropriate dataset in advance.

\section{Conclusion}\label{sec:conclusion}
This work introduces an approach for mixing constitutive models in finite element analysis and tests its performance by combining a physics-based constitutive model with a data-driven surrogate counterpart.
The method dynamically determines where to use each model based on the surrogate uncertainty, creating a spatially varying mixture that preserves accuracy while reducing computational cost.
This is achieved by solving a phase field based on the surrogate model uncertainty that determines which constitutive model is used.
Then, given the mixture of constitutive models, the mechanical problem of interest is solved.
By solving the problem in a staggered way a time-step dependency occurs, which can be resolved by performing multiple staggered iterations per load step.
Running several iterations brings additional costs but is necessary for the robustness of the method.

The use of a phase field was motivated by the hypothesis that the spatial smoothness would improve the numerical stability of the mechanical problem.
To test this, we studied the effect of decreasing the size of the interface.
As there is a limit to the minimum interface width for the phase field itself, we compared the phase-field approach to an alternative local approach with no spatial smoothening.
There, we used the surrogate model if its uncertainty is below a threshold, and the high-fidelity model otherwise, creating a sharp transition between the models.
We demonstrated on a notched plate problem that while the local approach can work and reduce the high-fidelity evaluations further, it is less robust than using a phase field.
For this reason, the spatial smoothness introduced by the phase field is beneficial.

The results demonstrate that the mixing approach significantly reduces the number of high-fidelity model evaluations compared to the full model.
The quality of the surrogate model directly impacts performance --- better-trained surrogate models with more data allow for increased computational savings.
If a surrogate model can be curated such that it accurately captures the behavior during the full simulation, then the uncertainty remains low, and no additional high-fidelity model evaluations are required.
Then, this framework adds the potential extra cost of requiring uncertainty from the surrogate.
However, creating an all-encompassing surrogate model is very challenging, and it is generally difficult to know in advance whether the surrogate model is sufficient.
Therefore, this mixture of models uses the surrogate where possible while still enabling an accurate simulation.
If the surrogate is poorly trained and can only capture a small part of the behavior, the additional Newton-Raphson iterations to obtain a consistent solution can result in more high-fidelity model evaluations than simply using the high-fidelity model without this approach.
It is, however, unlikely that this would be the case in practice.
This work thus provides a framework for accelerating multiscale simulations where micromodel evaluations are computationally expensive.
The approach works independently of the specific surrogate model and can, therefore, be combined with the latest advances in developing surrogate models to enable multiscale simulations.

\section*{Acknowledgements}\label{sec:acknowledge}
This work is supported by the TU Delft AI Labs programme.
We gratefully acknowledge the TU Delft AI Lab Travel grant for enabling this work.

% End of main text

% If you have bibdatabase file and want bibtex to generate the
% bibitems, please use
%
\bibliographystyle{elsarticle-num}
\bibliography{reflibproposal}

\appendix
\newpage

\section{Gaussian Process regression}\label{sec:Gaussianprocess}
A comprehensive treatment of GPs can be found in~\citep{williams1995gaussian}.
A Gaussian Process is a collection of random variables, any finite number of which follow a joint Gaussian distribution.
GPs are specified by a mean $m(\mathbf{x})$ and a covariance function $k(\mathbf{x}, \mathbf{x}')$, written as:
\begin{equation}
    f(\mathbf{x}) \sim \mathrm{GP}(m(\mathbf{x}), k(\mathbf{x}, \mathbf{x}')).
\end{equation}
For simplicity, it is common to assume a zero mean function, $m(\mathbf{x})=0$.
The covariance function $k(\mathbf{x}, \mathbf{x}')$ specifies the correlation between outputs at different inputs, often capturing how similar inputs should yield similar outputs.
In this work, we use a squared exponential (radial basis) kernel function:
\begin{equation}
    k(\mathbf{x}, \mathbf{x}') = \sigma^2_f \exp\left(-\frac{\|\mathbf{x} - \mathbf{x}'\|^2}{2\ell^2}\right),
\end{equation}
where $\sigma^2_f$ is the signal variance, and $\ell$ is the characteristic length scale.

Given a training dataset $\mathcal{D}=\{(\mathbf{x}_i, y_i) | i=1,...,n\}$ with inputs $\mathbf{X} = [\mathbf{x}_1,...,\mathbf{x}_n]^T$ and outputs $\mathbf{y} = [y_1,...,y_n]^T$, we aim to predict the output $f(\mathbf{x}_*)$ at a new test point $\mathbf{x}_*$.
Assuming noisy observations $y_i = f(\mathbf{x}_i) + \epsilon_i$ where $\epsilon_i \sim \mathcal{N}(0, \sigma^2_n)$, the joint distribution of the training outputs and the test output becomes:
\begin{equation}
    \begin{bmatrix} \mathbf{y} \\ f(\mathbf{x}_*) \end{bmatrix} \sim \mathcal{N}\left( \mathbf{0}, \begin{bmatrix} K(\mathbf{X}, \mathbf{X}) + \sigma^2_n\mathbf{I} & K(\mathbf{X}, \mathbf{x}_*) \\ K(\mathbf{x}_*, \mathbf{X}) & K(\mathbf{x}_*, \mathbf{x}_*) \end{bmatrix} \right),
\end{equation}
where $K(\mathbf{X}, \mathbf{X})$ denotes the $n \times n$ matrix of covariances evaluated at all pairs of training points, $K(\mathbf{X}, \mathbf{x}_*)$ is the $n \times 1$ vector of covariances between training points and the test point, and $K(\mathbf{x}_*, \mathbf{x}_*)$ is the covariance of the test point with itself.
The posterior distribution of $f(\mathbf{x}_*)$ given the observations becomes:
\begin{align}
    p(f(\mathbf{x}_*) | \mathbf{X}, \mathbf{y}, \mathbf{x}_*) &= \mathcal{N}(\mu_*, \sigma^2_*) \\
    \mu_* &= K(\mathbf{x}_*, \mathbf{X})[K(\mathbf{X}, \mathbf{X}) + \sigma^2_n\mathbf{I}]^{-1}\mathbf{y} \\
    \sigma^2_* &= K(\mathbf{x}_*, \mathbf{x}_*) - K(\mathbf{x}_*, \mathbf{X})[K(\mathbf{X}, \mathbf{X}) + \sigma^2_n\mathbf{I}]^{-1}K(\mathbf{X}, \mathbf{x}_*)
\end{align}
The predictive mean $\mu_*$ serves as our point estimate, while the predictive variance $\sigma^2_*$ quantifies our uncertainty about the estimate.
This uncertainty naturally increases as we move away from the training data.

There are several methods for optimizing the hyperparameters of a GP.
We use the maximum likelihood estimation approach, which finds the hyperparameters that maximize the marginal likelihood of the observed data.
The hyperparameters $\boldsymbol{\theta}$ (i.e. $\boldsymbol{\theta} = \{\sigma_f, \ell, \sigma_n\}$) are optimized by maximizing the log marginal likelihood:
\begin{equation}
    \log p(\mathbf{y}|\mathbf{X}, \boldsymbol{\theta}) = -\frac{1}{2}\mathbf{y}^T[K(\mathbf{X}, \mathbf{X}) + \sigma^2_n\mathbf{I}]^{-1}\mathbf{y} - \frac{1}{2}\log|K(\mathbf{X}, \mathbf{X}) + \sigma^2_n\mathbf{I}| - \frac{n}{2}\log(2\pi)
\end{equation}
We perform this optimization using the L-BFGS-B algorithm, starting from 20 different initial configurations to avoid local optima.

\end{document}